\documentclass[11pt]{article}
\usepackage{amsmath}
\usepackage{amssymb}
\usepackage{amscd}
\usepackage{latexsym}
\usepackage{theorem}
\usepackage{doublespace}
\usepackage{epsfig}
\usepackage{psfrag}

\newtheorem{theorem}{Theorem}[section]
\newtheorem{lemma}[theorem]{Lemma}
\newtheorem{proposition}[theorem]{Proposition}
\newtheorem{corollary}[theorem]{Corollary}

{\theorembodyfont{\rmfamily}
\theoremstyle{plain}
\newtheorem{definition}[theorem]{Definition}
\newtheorem{example}[theorem]{Example}

\newtheorem{remark}[theorem]{Remark}
}

\input xy
\xyoption{all}

\newcommand{\Hom}{\operatorname{Hom}}
\renewcommand{\Im}{\operatorname{Im}}

\newcommand{\id}{\operatorname{id}}

\newcommand{\Pd}{\displaystyle\prod}
\newcommand{\Sd}{\displaystyle\sum}

\newcommand{\C}{{\mathbb{C}}}
\newcommand{\Z}{{\mathbb{Z}}}
\newcommand{\Q}{{\mathbb{Q}}}
\newcommand{\R}{{\mathbb{R}}}
\renewcommand{\P}{{\mathbb{P}}}
\renewcommand{\H}{{\mathbb{H}}}

\renewcommand{\a}{\alpha}
\renewcommand{\b}{\beta}
\newcommand{\g}{\gamma}

\newcommand{\K}{\kappa}

\newcommand{\Kh}{\hat\K}
\newcommand{\Hn}{\H^n}
\newcommand{\gd}{\mathfrak{g}^*}
\newcommand{\td}{\mathfrak{t}^*}
\newcommand{\rgt}{r_T^G}

\newcommand{\half}{\frac{1}{2}}

\newcommand{\mur}{\mu_{\R}}
\newcommand{\muc}{\mu_{\C}}
\newcommand{\muh}{\mu_{G}}
\newcommand{\mubh}{\mu_{T}}

\newcommand{\Cn}{\C^n}

\renewcommand{\cot}{T^*\Cn}

\renewcommand{\mod}{{/\!\!/}}
\newcommand{\mmod}{{/\!\!/\!\!/\!\!/}}
\newcommand{\bigmod}{{\Big/\!\!\!\Big/}}
\newcommand{\bigmmod}{{\Big/\!\!\!\Big/\!\!\!\Big/\!\!\!\Big/}}

\newcommand{\hookto}{{\hookrightarrow}}
\renewcommand{\iff}{\Leftrightarrow}
\newcommand{\impl}{\Rightarrow}
\newcommand{\subs}{\subseteq}

\newcommand{\ra}{\rangle}
\newcommand{\la}{\langle}
\newcommand{\IM}{\int_M}
\newcommand{\IN}{\int_N}
\newcommand{\IG}{\int_G}
\newcommand{\IF}{\int_F}
\newcommand{\intxt}{\int_{X\mmod T}}
\newcommand{\intxg}{\int_{X\mmod G}}

\newcommand{\I}{\mathcal{I}}
\newcommand{\J}{\mathcal{J}}
\newcommand{\imkt}{\left(\Im\K_T\right)}
\newcommand{\so}{S^1}
\newcommand{\hso}{H_{\so}^*}
\newcommand{\hhso}{\widehat{H}_{\so}^*}
\newcommand{\hxg}{H_{\so}^*(X\mmod G)}
\newcommand{\hxt}{H_{\so}^*(X\mmod T)}
\newcommand{\hhxg}{\widehat{H}_{\so}^*(X\mmod G)}
\newcommand{\hhxt}{\widehat{H}_{\so}^*(X\mmod T)}

\newcommand{\hhsg}{\widehat{H}_{\so\times G}^*(X)}
\newcommand{\hp}{H_{S^1}^*(pt)}
\newcommand{\hhp}{{\mathbb K}}
\newcommand{\hm}{\hso(M)}
\newcommand{\hhn}{\hhso(N)}
\newcommand{\hhm}{\hhso(M)}
\newcommand{\M}{{\mathfrak{M}}}

\newcommand{\NN}{\mathfrak{N}}
\newcommand{\dso}{\Delta_*(1)}

\newcommand{\hs}{\hspace{3pt}}

\newcommand{\qed}{\hfill \mbox{$\Box$}\medskip\newline}
\newenvironment{proof}{\noindent {\bf Proof:}}{\qed \par}
\newenvironment{proofintegration}{\noindent {\bf Proof of \ref{integration}:}}{\qed \par}
\newenvironment{proofordinary}{\noindent {\bf Proof of \ref{ordinary}:}}{\qed \par}
\newenvironment{proofmain}{\noindent {\bf Proof of \ref{main}:}}{\qed \par}

\begin{document}
\begin{spacing}{1.1}

\noindent
{\Large \bf Abelianization for hyperk\"ahler quotients}
\bigskip\\
{\bf Tam\'as Hausel} \\
Department of Mathematics, University of Texas,
Austin, TX 78712\smallskip \\
{\bf Nicholas Proudfoot } \\
Department of Mathematics, University of California,
Berkeley, CA 94720
\bigskip
{\small
\begin{quote}
\noindent {\em Abstract.}
We study an integration theory in circle equivariant cohomology
in order to prove a theorem relating the cohomology ring of a hyperk\"ahler
quotient to the cohomology ring of the quotient by a maximal abelian
subgroup, analogous to a theorem of Martin for symplectic quotients.
We discuss applications of this theorem to quiver varieties, and
compute as an example the ordinary and equivariant cohomology rings of 
a hyperpolygon space.
\end{quote}
}
\bigskip

Let $X$ be a symplectic manifold equipped with a hamiltonian
action of a compact Lie group $G$.  Let $T\subs G$
be a maximal torus, let $\Delta\subset\td$ be the set of roots of $G$,
and let $W=N(T)/T$ be the Weyl group.  
If the symplectic quotients
$X\mod G$ and $X\mod T$ are both compact, Martin's theorem \cite[Theorem A]{M}
relates the cohomology\footnote{In this paper cohomology means cohomology with rational 
coefficients.} of $X\mod G$ to the cohomology of $X\mod T$.
Specifically, it says that $$H^*(X\mod G)\cong 
\frac{H^*(X\mod T)^W}{ann(e_0)},$$
where $$e_0 = \prod_{\a\in\Delta}\a\in\left(\operatorname{Sym}\td\right)^W
\cong H^*_T(pt)^W,$$ which acts naturally on
$H^*(X\mod T)^W \cong H^*_T(\mu_T^{-1}(0))^W$. 
In the case where $X$ is a complex vector space and $G$ acts linearly
on $X$, a similar result was obtained by Ellingsrud and Str\o mme \cite{ES}
using different techniques.

Our goal is to state and prove an analogue of this theorem for
hyperk\"ahler quotients \cite{HKLR}.  There are two main obstacles to this goal.
First, hyperk\"ahler quotients are rarely compact.
The assumption of compactness in Martin's theorem is crucial
because his proof involves integration. 
Generalizing an idea of \cite{MNS} and \cite{P}, 
our answer to this problem is to work
with equivariant cohomology
of {\em circle compact} manifolds, by which we mean oriented manifolds
with an action of $\so$ such that the fixed point 
set is oriented and compact.   
By the localization theorem of Atiyah-Bott \cite{AB} and Berline-Vergne 
\cite{BV},
integration in rationalized $\so$-equivariant cohomology
of circle compact manifolds can be defined in terms of integration
on their fixed point sets.  Section~\ref{thepush} is devoted to
making this statement precise by defining a well-behaved 
push forward in the  rationalized $\so$-equivariant cohomology 
of circle compact manifolds.

The second obstacle is that Martin's result uses surjectivity \cite{Ki}
of the Kirwan map from $H^*_G(X)$ to $H^*(X\mod G)$.  The analogous
map for hyperk\"ahler quotients is surjective only conjecturally.
Our approach is to assume that the rationalized Kirwan map is
surjective, which is equivalent to saying that the cokernel
of the non-rationalized Kirwan map
$$\K_G:H^*_{\so\times G}(X)\to\hxg$$
is torsion as a module over $\hp$.  This is a weaker assumption
than surjectivity of $K_G$; in particular, we show in Section
\ref{quiver} that this assumption holds for quiver varieties,
as a consequence of the work of Nakajima.

Under this assumption,
Theorem~\ref{main} computes the rationalized equivariant cohomology
of $X\mmod G$ in terms of that of $X\mmod T$. 
We show that
$$\hhso(X\mmod G)\cong\frac{\hhso(X\mmod T)^W}{ann(e)},$$
where
$$e = \Pd_{\a\in\Delta}\a(x-\a) \in (\operatorname{Sym}\td)^W\otimes\Q[x]
\cong H^*_{\so\times T}(pt)^W.$$
Theorem~\ref{ordinary}
describes the image of the non-rationalized Kirwan map
in a similar way:
$$\hso(X\mmod G)\supseteq\Im(\K_G)\cong\frac{\imkt^W}{ann(e)},$$
where $\K_T:H^*_{\so\times T}(X)\to \hso(X\mmod T)$ is the Kirwan
map for the abelian quotient.
In many situations, such as when $X=\cot$,
$\K_T$ is known to be surjective.

In Section~\ref{quiver} we show that all of the hypotheses
of Theorems~\ref{main} and \ref{ordinary} are satisfied for
Nakajima's quiver varieties. This way we can reduce questions
about the (rationalized) equivariant cohomology of quiver varieties
to questions about the (rationalized) equivariant cohomology of 
toric hyperk\"ahler varieties (also called hypertoric varieties in \cite{HP1}). 
The cohomology rings of toric hyperk\"ahler 
varieties are well understood, as in \cite{BD}, \cite{HP1}, \cite{HS} and 
\cite{K1}.  When the hyperk\"ahler Kirwan map is known to 
be surjective, for example in the case of the Hilbert scheme of points on an ALE space, 
Theorem \ref{ordinary} gives an explicit description
of the cohomology ring of the quiver variety.  
Such cases are discussed in Remarks \ref{generation} and \ref{horn}.

We conclude in Section~\ref{polygon} by demonstrating how the ideas of the present paper work
in the case of a particular quiver variety, the so-called hyperpolygon space. 
We show that the hyperk\"ahler Kirwan map is surjective, and therefore our machinery
reproduces, by different means, the results of \cite[\S 7]{K2} and \cite[\S 3]{HP2}.

\paragraph{\bf Acknowledgment.} We would like to acknowledge useful conversations
with Hiraku Nakajima and Michael Thaddeus. In particular an example of Thaddeus
is used in Example~\ref{thaddeus}. Financial support
was provided in part by NSF grants DMS-0072675 and DMS-0305505. 

\begin{section}{Integration}\label{thepush}
The localization theorem of Atiyah-Bott \cite{AB} and Berline-Vergne \cite{BV} says that
given a manifold $M$ with a circle action, the restriction map from
the circle equivariant cohomology of $M$ to the circle equivariant
cohomology of the fixed point set $F$ is an isomorphism modulo torsion.
In particular, integrals on a compact $M$ can be computed in terms of
integrals on $F$.
If $F$ is compact, it is possible to use the Atiyah-Bott-Berline-Vergne formula
to {\em define} integrals on $M$.

We will work in the category of {\it circle compact}
manifolds, by which we mean oriented $\so$-manifolds with compact
and oriented fixed point sets.  
Maps between circle compact
manifolds are required to be equivariant.

\begin{definition}
Let $\hhp = \Q(x)$, the rational function field of $\hp\cong \Q[x]$.
For a circle compact manifold $M$, let $\hhm = \hm\otimes\hhp$,
where the tensor product is taken over the ring $\hp$. We call $\hhm$
the {\em rationalized} $\so$-equivariant cohomology of $M$. 
\end{definition}

An immediate consequence of \cite{AB} is that restriction gives an isomorphism 
\begin{equation}\label{isomorphic}\hhm\cong \hhso(F)\cong H^*(F)\otimes_{\Q}\hhp,\end{equation} 
where $F=M^{S^1}$ denotes 
the compact fixed point set of $M$.  In particular
$\hhm$ is a finite dimensional vector space over $\hhp$, 
and trivial if and only if $F$ is empty. 

Let $i:N\hookto M$ be a closed embedding.  
There is a standard notion of proper pushforward
$$i_*:\hso(N)\to\hm$$ given by the formula $i_*=r\circ\Phi$,
where $r:\hso(M,M\setminus N)\to\hso(M)$ is the restriction map,
and $\Phi:\hso(N)\to\hso(M,M\setminus N)$ is the Thom isomorphism.
We will also denote the induced map $\hhso(N)\to\hhm$ by $i_*$.
Geometrically, $i_*$ can be understood as the inclusion
of cycles in Borel-Moore homology.

This map satisfies two important formal properties
\cite{AB}:
\begin{equation}\label{functor}
\text{Functoriality:   } (i\circ j)_* = i_*\circ j_*
\end{equation}
\begin{equation}\label{module}
\text{Module homomorphism:   } i_*(\g\cdot i^*\a)
=i_*\g\cdot\a 
\text{   for all   }\a\in\hhm, \g\in\hhn.
\end{equation}
We will denote the Euler
class $i^*i_*(1)\in\hhso(N)$ by $e(N)$.
If a class $\g\in\hhn$ is in the image of $i^*$,
then property (\ref{module}) tells us that
$i^*i_*\g = e(N)\g$.  Since the pushforward construction is 
local in a neighborhood
of $N$ in $M$, we may assume that $i^*$ is surjective, hence this identity
holds for all $\g\in\hhn$.

Let $F= M^{\so}$ be the fixed point set of $M$. 
Since $M$ and $F$ are each oriented, so is the normal bundle
to $F$ inside of $M$.  The following result is 
standard, see e.g. \cite{Ki}. 

\begin{lemma}
The Euler class $e(F)\in\hhso(F)$ of the normal bundle to $F$ in $M$ 
is invertible.
\end{lemma}

\begin{proof}
Since $\so$ acts trivially on $F$, $\hhso(F)\cong H^*(F)\otimes_{\Q}\hhp$.
We have $e(F) = 1\otimes ax^k + nil$, where
$k= \operatorname{codim}(F)$, $a$ is the product of the weights of the $S^1$
action on any fiber of the normal bundle, and $nil$ consists of terms
of positive degree in $H^*(F)$.  Since $F$ is the fixed
point set, $\so$ acts freely on the complement of the zero section
of the normal bundle, therefore $a\neq 0$.  Since $ax^k$ is invertible and
$nil$ is nilpotent, we are done.
\end{proof}

\vspace{-\baselineskip}
\begin{definition}
For $\a\in\hhm$, let $$\IM\a=\IF\frac{\a|_F}{e(F)}\in\hhp.$$
\end{definition}

Note that this definition does not depend on our choice
of orientation of $F$.  Indeed, reversing the orientation of $F$
has the effect of negating $e(F)$, {\em and} introducing a second factor
of $-1$ coming from the change in fundamental class.  These two effects
cancel.

For this definition to be satisfactory, we must be able to prove
the following lemma, which is standard in the setting of
ordinary cohomology of compact manifolds.

\begin{lemma}\label{Npush}
Let $i:N\hookto M$ be a closed immersion.
Then for any $\a\in\hhm,\g\in\hhso(N)$, we have
$\IM\a\cdot i_*\g = \IN i^*\a\cdot\g$.
\end{lemma}

\begin{proof}
Let $G=N^{\so}$, let $j:G\to F$ denote the restriction of $i$ to $G$,
and let $\phi:F\to M$ and $\psi:G\to N$ denote the inclusions of $F$ and $G$
into $M$ and $N$, respectively.
$$\begin{CD}
N        @>i>> M\\
@A\psi AA @AA\phi A\\
G        @>j>> F
\end{CD}$$
Then $$\IM\a\cdot i_*\g = \IF\frac{\phi^*\a\cdot\phi^*i_*\g}{e(F)},$$
and $$\IN i^*\a\cdot\g = \IG\frac{\psi^*i^*\a\cdot\psi^*\g}{e(G)}
=\IG\frac{j^*\phi^*\a\cdot\psi^*\g}{e(G)}
=\IF\phi^*\a\cdot j_*\left(\frac{\psi^*\g}{e(G)}\right),$$
where the last equality is simply the integration formula applied
to the map $j:G\to F$ of compact manifolds \cite{AB}.
Hence it will be sufficient to prove that
$$\phi^*i_*\g=e(F)\cdot j_*\left(\frac{\psi^*\g}{e(G)}\right)\in\hhso(F).$$
To do this, we will
show that the difference of the two classes lies in
the kernel of $\phi_*$, which we know is trivial because the composition
$\phi^*\phi_*$ is given by multiplication by the invertible 
class $e(F)\in\hhso(F)$.
On the left hand side we get
$$\phi_*\phi^*i_*\g = \phi_*(1)\cdot i_*\g\hspace{15pt}
\text{by    }(\ref{module}),
$$
and on the right hand side we get
\begin{eqnarray*}
\phi_*\left(e(F)\cdot j_*\left(\frac{\psi^*\g}{e(G)}\right)\right)
&=& \phi_*\left(\phi^*\phi_*(1)
\cdot j_*\left(\frac{\psi^*\g}{e(G)}\right)\right)\\
&=& \phi_*(1)\cdot \phi_*j_*\left(\frac{\psi^*\g}{e(G)}\right)
\hspace{1in}\text{    by    }(\ref{module})\\
&=& \phi_*(1)\cdot i_*\psi_*\left(\frac{\psi^*\g}{e(G)}\right)
\hspace{1in}\text{    by    }(\ref{functor}).
\end{eqnarray*}
It thus remains only to show that
$\g = \psi_*\left(\frac{\psi^*\g}{e(G)}\right)$.
This is seen by applying $\psi^*$ to both sides,
which is an isomorphism (working over the field $\hhp$) by \cite{AB}.
\end{proof}

For $\a_1,\a_2\in\hhm$, consider the
symmetric, bilinear, $\hhp$-valued pairing
$$\la\a_1,\a_2\ra_M=\IM\a_1\a_2.$$

\begin{lemma}[Poincar\'e Duality]\label{perfect}
This pairing is nondegenerate.
\end{lemma}

\begin{proof}
Suppose that $\a\in\hhm$ is nonzero, and therefore $\phi^*\a\neq 0$.  
Since $F$ is
compact, there must exist a class $\g\in\hhso(F)$ such
that $0\neq\IF\phi^*\a\cdot\g=\IM\a\cdot \phi_*\g=\la\a,\phi_*\g\ra_M$.
\end{proof}

\vspace{-\baselineskip}
\begin{definition}\label{push}
For an arbitrary equivariant map $f:N\to M$, we may now define
the pushforward $$f_*:\hhn\to\hhm$$ to
be the adjoint of $f^*$ with respect to the pairings
$\la\cdot,\cdot\ra_N$ and $\la\cdot,\cdot\ra_M$.  
This is well defined because, according to (\ref{isomorphic}), 
$\hhm$ and $\hhn$ are finite dimensional
vector spaces over the field $\hhp$. Lemma~\ref{Npush} tells us that this definition generalizes
the definition for closed immersions.  Furthermore,
properties (\ref{functor}) and (\ref{module}) for pushforwards along 
arbitrary maps
are immediate corollaries of the definition.
If $f$ is a projection, then $f_*$ will be given by 
integration along the fibers.
Using the fact that 
every map factors through its graph as a closed immersion 
and a projection, we always have a geometric interpretation of the pushforward.
\end{definition}

As an application, let us consider the manifold $M\times M$,
along with the two projections $\pi_1$ and $\pi_2$,
and the diagonal map $\Delta:M\to M\times M$.
Suppose that we can write $$\Delta_*(1)=\sum\pi_1^*\a_i\cdot\pi_2^*\b_i$$
for a finite collection of classes $a_i,b_i\in\hhm$.
The following Proposition will be used in Section~\ref{quiver}.

\begin{proposition}\label{basis}
The set $\{b_i\}$ is an additive basis for $\hhm$.
\end{proposition}

\begin{proof}
For any $\a\in\hhm$, we have
\begin{eqnarray*}
\a &=& \id_*\id^*\a\\
&=& (\pi_2\circ\Delta)_*(\pi_1\circ\Delta)^*\a\\
&=& \pi_{2*}\big(\Delta_*\left(1\cdot\Delta^*\pi_1^*\a\right)\big)\\
&=& \pi_{2*}\big(\pi_1^*\a\cdot\Delta_*(1)\big)\\
&=& \pi_{2*}\left(\sum\pi_1^*(a_i\a)\cdot\pi_2^*b_i\right)\\
&=& \sum\pi_{2*}\pi_1^*(a_i\a) \cdot b_i\\
&=& \sum\la a_i,\a\ra\cdot b_i,
\end{eqnarray*}
hence $\a$ is in the span of $\{b_i\}$.
\end{proof}
\end{section}

\begin{section}{An analogue of Martin's theorem}\label{martinsection}
Let $X$ be a hyperk\"ahler manifold with a circle action, and
suppose that a compact Lie group $G$ acts hyperhamiltonianly on $X$.
We will assume that the circle action preserves a given complex structure
$I$.
Having chosen a particular complex structure on $X$, we may
write the hyperk\"ahler moment map in the form
$$\muh=\mur\oplus\muc:X\to\gd\oplus\gd_{\C},$$
where $\muc$ is holomorphic with respect to $I$ \cite{HP1}.
We require that the action of $G$ commute with the action of $S^1$,
that $\mur$ is $\so$-invariant, and that $\muc$ is $\so$-equivariant
with respect to the action of $S^1$ on $\gd_{\C}$ by complex multiplication.
We do {\em not} ask the action of $S^1$ on $X$ to preserve the hyperk\"ahler
structure.

Let $T\subs G$ be a maximal torus, and let $pr:\gd\to\td$
be the natural projection.  Then $T$ acts on $X$ with hyperk\"ahler
moment map
$$\mubh=pr\circ\mur\oplus pr_{\C}\circ\muc:
X\to\td\oplus\td_{\C}.$$
Let $\xi\in\gd$ be a central element such that $(\xi,0)$ is a regular
value of $\muh$ and $(pr(\xi),0)$ is a regular value of $\mubh$.
Let $$X\mmod G = \muh^{-1}(\xi,0)/G\hs\hs\hs\text{   and   }\hs\hs\hs
X\mmod T = \mubh^{-1}(pr(\xi),0)/T$$
be the hyperk\"ahler quotients of $X$ by $G$ and $T$, respectively.
Because $\muh$ and $\mubh$ are circle equivariant,
the action of $S^1$ on $X$ descends to actions on the hyperk\"ahler quotients.
Note that $X\mmod T$ also inherits an action of the Weyl group $W$ of $G$.

\begin{example}\label{standard}
Suppose that $G$ acts linearly on $\Cn$ with moment map
$\mu:\Cn\to\gd$, and let $X$ be the hyperk\"ahler manifold
$\cot\cong\Hn$.  The action of $G$ on $\Cn$ induces an action
of $G$ on $X$ with hyperk\"ahler moment map
$$\mur(z,w) = \mu(z)-\mu(w)\text{  and  }
\muc(z,w)(v) = w(\hat{v}_z),$$
where $w\in T_z^*\Cn\cong\Cn$, $v\in\gd_{\C}$, and
$\hat{v}_z$ the element of $T_z\Cn$ induced by $v$ \cite{HP1}.
The action of $G$ commutes with the action of $S^1$ on $X$
given by scalar multiplication on each fiber, and the
hyperk\"ahler moment map is equivariant.
The quotient $X\mmod G$
is a partial compactification of the cotangent bundle
$T^*(X\mod G)$, and is circle compact if $\mu$
is proper \cite[1.3]{HP1}.
\end{example}

Consider the Kirwan maps
$$\K_G:H_{\so\times G}^*(X)\to \hso(X\mmod G)\hs\hs\hs\text{   and   }\hs\hs\hs
\K_T:H_{\so\times T}^*(X)\to \hso(X\mmod T),$$
induced by the inclusions of $\muh^{-1}(\xi,0)$ and $\mubh^{-1}(pr(\xi),0)$ into $X$,
along with their rationalizations
$$\hat\K_G:\widehat H_{\so\times G}^*(X)\to \hhso(X\mmod G)\hs\hs\hs\text{   and   }\hs\hs\hs
\hat\K_T:\widehat H_{\so\times T}^*(X)\to \hhso(X\mmod T).$$
Let $$\rgt:\widehat{H}^*_{\so\times G}(X)\to\widehat{H}^*_{\so\times T}(X)^W$$
be the standard isomorphism.

Let $\Delta=\Delta^+\sqcup\Delta^-\subset\td$ be the set of roots of $G$.
Let $$e = \Pd_{\a\in\Delta}\a(x-\a) \in (\operatorname{Sym}\td)^W\otimes\Q[x]
\cong H_{\so\times G}(pt)\subs\widehat{H}_{\so\times G}(pt),$$
and $$e' = \Pd_{\a\in\Delta^-}\a\cdot\Pd_{\a\in\Delta}(x-\a)
\in\operatorname{Sym}\td\otimes\Q[x]\cong H_{\so\times T}(pt)
\subs\widehat{H}_{\so\times T}(pt).$$
The following two theorems are analogues of
Theorems B and A of \cite{M}, adapted to circle compact
hyperk\"ahler quotients.  Our proofs follow closely those of Martin.

\begin{theorem}\label{integration}
Suppose that $X\mmod G$ and $X\mmod T$ are both circle compact.
If $\g\in\hhsg$, then
$$\intxg \Kh_G(\g) = \frac{1}{|W|}\intxt \Kh_T\circ\rgt(\g)\cdot e.$$
\end{theorem} 

\begin{theorem}\label{main}
Suppose that $X\mmod G$ and $X\mmod T$ are both circle compact,
and that the rationalized Kirwan map $\Kh_G$ 
surjective.  Then
$$\hhso(X\mmod G)\cong\frac{\hhso(X\mmod T)^W}{ann(e)}
\cong\left(\frac{\hhso(X\mmod T)}{ann(e')}\right)^W.$$
\end{theorem}

\begin{proofintegration}
Consider the following pair of maps:
$$\begin{CD}
\muh^{-1}(\xi,0)/T  @>i>> \mubh^{-1}(pr(\xi),0)/T\cong X\mmod T\\
@V\pi VV \\
\muh^{-1}(\xi,0)/G\cong X\mmod G.
\end{CD}$$
Each of these spaces is a complex $S^1$-manifold with a compact, complex
fixed point set, and therefore satisfies the hypotheses of Section
\ref{thepush}.
Let $$b = \Pd_{\a\in\Delta^+}\a\in H_{\so\times T}(pt)$$ be the product of the
positive roots of $G$, which we will think of as an
element of $\hhxt$.  Martin shows that $\pi_*i^*b=|W|$, 
and that $i^*\circ\Kh_T\circ\rgt = \pi^*\Kh_G$ \cite{M}, hence
we have
\begin{eqnarray*}
\intxg\Kh_G(\g) &=& \frac{1}{|W|}\intxg\Kh_G(\g)\cdot\pi_*i^*b\\
&=& \frac{1}{|W|}\int_{\muh^{-1}(\xi,0)/T}\pi^*\Kh_G(\g)\cdot i^*b\hspace{20pt}
\text{by Definition~\ref{push}}\\
&=& \frac{1}{|W|}\int_{\muh^{-1}(\xi,0)/T}i^*\circ\Kh_T\circ\rgt(\g)\cdot i^*b\\
&=& \frac{1}{|W|}\intxt \Kh_T\circ\rgt(\g)\cdot b \cdot i_*(1)\hspace{20pt}
\text{by Lemma~\ref{Npush}}.
\end{eqnarray*}
It remains to compute $i_*(1)\in\hhxt$.
For $\a\in\Delta$, let $$L_{\a} = \mubh^{-1}((pr(\xi),0)\times_T\C_{\a}$$
be the line bundle on $X\mmod T$ with $S^1$-equivariant
Euler class $\a$.
Similarly, let $L_x$ be the (topologically trivial)
line bundle with $S^1$-equivariant
Euler class $x$.
Following the idea of \cite[1.2.1]{M}, we observe
that the restriction of $\muh - (\xi,0)$ to 
$\mubh^{-1}(pr(\xi),0)$ defines
an $\so\times T$-equivariant map $$s:\mubh^{-1}(pr(\xi),0)
\to V\oplus V_{\C},$$ where $V = pr^{-1}(0)$ 
and $V_{\C}=pr_{\C}^{-1}(0)$.
This descends to an $S^1$-equivariant section of the
bundle $E = \mubh^{-1}(pr(\xi),0)\times_T \left(V\oplus V_{\C}\right)$
with zero locus $\muh^{-1}(\xi,0)/T$.  The fact that $(\xi,0)$
is a regular value implies that this section is generic,
hence the equivariant Euler class $e(E)\in\hhxt$ is equal to $i_*(1)$.

The vector space $V$ is isomorphic as a $T$-representation to 
$\bigoplus_{\a\in\Delta^-}\C_{\a}$,
with $S^1$ acting trivially.
Similarly, $V_{\C}$ is isomorphic to $V\otimes\C\cong V\oplus V^*$,
with $S^1$ acting diagonally by scalars.
Hence 
\begin{eqnarray*}
E &\cong& \bigoplus_{\a\in\Delta^-}L_{\a}
\oplus\bigoplus_{\a\in\Delta^-}\left(L_x\otimes L_{\a}\right)
\oplus \left(L_x\otimes L_{-\a}\right)\\
&\cong& \bigoplus_{\a\in\Delta-}L_{\a}\oplus
\bigoplus_{\a\in\Delta}L_x\otimes L_{-\a},
\end{eqnarray*}
and therefore $$i_*(1)=e(E) = \prod_{\a\in\Delta^-}\a\cdot\prod_{\a\in\Delta}(x-\a)=e'.$$
Multiplying by $b$ we obtain $e$, and the theorem is proved.
\end{proofintegration}

\begin{proofmain}
Observe that the restriction
of $\pi^*$ to the Weyl-invariant part 
$\hhso\left(\muh^{-1}(\xi,0)/T\right)^W$ is given by the composition of
isomorphisms
$$\hhso\left(\muh^{-1}(\xi,0)/T\right)^W
\cong\widehat{H}^*_{\so\times T}\left(\muh^{-1}(\xi,0)\right)^W
\cong\widehat{H}^*_{\so\times G}\left(\muh^{-1}(\xi,0)\right)
\cong\hhxg,$$
hence we may define
$$i^*_W:= (\pi^*)^{-1}\circ i^*:
\hhxt^W\to\hhso\left(\muh^{-1}(\xi,0)/T\right)^W.$$
Furthermore, we have
$\Kh_G = i^*_W\circ\Kh_T\circ\rgt$, hence
$i^*_W$ is surjective.
As in \cite[\S 3]{M},
\begin{eqnarray*}
i^*_W(a)=0 &\iff& \forall c\in\hhxt^W, \intxg i^*_W(c)\cdot i^*_W(a)=0\hspace{15pt}
\text{by \ref{perfect} and surjectivity of $i^*_W$}\\
&\iff& \forall c\in\hhxt^W, \intxt c\cdot a\cdot e = 0\hspace{15pt}
\text{by Theorem~\ref{integration}}\\
&\iff& \forall d\in\hhxt, \intxt d\cdot a\cdot e=0\hspace{15pt}
\text{by using $W$ to average $d$}\\
&\iff& a\cdot e=0\hspace{15pt}\text{by Lemma~\ref{perfect}},
\end{eqnarray*}
hence $\ker i^*_W = ann(e)$.
By surjectivity of $i^*_W$,
\begin{eqnarray*}
\hhso(X\mmod G) &\cong&
\frac{\hhxt^W}{\ker i^*_W}\cong\frac{\hhxt^W}{ann(e)}.
\end{eqnarray*}
By a second application of Lemma~\ref{perfect}, for any $a\in\hhxt$, we have
\begin{eqnarray*}
i^*(a) = 0 &\impl& \forall f\in\hhso(\muh^{-1}(\xi,0)/T),
\int_{\muh^{-1}(\xi,0)/T}f\cdot i^*(a) = 0\\
&\impl& \forall c\in\hhxt,
\int_{\muh^{-1}(\xi,0)/T} i^*(c)\cdot i^*(a)=0\\
&\impl& \forall c\in\hhxt,
\intxt c\cdot a\cdot i_*(1)=0\hspace{15pt}\text{by Lemma~\ref{Npush}}\\
&\impl& a\cdot e'= a\cdot i_*(1) = 0\hspace{15pt}\text{by Lemma~\ref{perfect}},
\end{eqnarray*}
hence $\ker i^* \subs ann(e')$.
This gives us a natural surjection
$$\frac{\hhxt^W}{ann(e)}=\frac{\hhxt^W}{\ker i^*_W}
\cong \left(\frac{\hhxt}{\ker i^*}\right)^W
\to \left(\frac{\hhxt}{ann(e')}\right)^W,$$
which is also injective because $e'$ divides $e$.
This completes the proof of Theorem \ref{main}.
\end{proofmain}

For the non-rationalized version of Theorem~\ref{main}, we make the additional
assumption that $X\mmod G$ and $X\mmod T$ are {\em equivariantly formal}
$\so$-manifolds, i.e. that $\hxg$ and $\hxt$ are free modules
over $\hp$.  This is the case whenever the circle action
is hamiltonian and its moment map is proper and bounded below 
(see \cite{Ki} and \cite[4.7]{HP1}).

\begin{theorem}\label{ordinary}
Suppose that $X\mmod G$ and $X\mmod T$ are equivariantly formal, circle compact, and 
that the rationalized Kirwan map $\Kh_G$ is surjective.
Then
$$\hso(X\mmod G)\supseteq\Im(\K_G)\cong\frac{(\Im\K_T)^W}{ann(e)}
\cong\left(\frac{\Im\K_T}{ann(e')}\right)^W.$$
\end{theorem}

\begin{remark}\label{notsobad}
In the context of Example~\ref{standard} with $pr\circ\mu$ proper, 
$X\mmod G$ and $X\mmod T$ are both circle compact and equivariantly
formal, and $\K_T$ is always surjective \cite{HP1}.
Note that this applies throughout Sections~\ref{quiver} and~\ref{polygon}.
\end{remark}

\begin{proofordinary}
Consider the following exact commutative diagram
$$
\xymatrix{
0 \ar[r] & A \ar[r]\ar[d] & \hxt^W \ar[r]^{i_W^*} \ar[d]
& \hxg \ar[d]\\
0 \ar[r] & \widehat{A} \ar[r] & \hhxt^W \ar[r]^{i_W^*} & \hhxg.}
$$
Equivariant formality implies that the downward maps in the above 
diagram are inclusions, hence the map on top labeled $i^*_W$
is simply the restriction of the map on the bottom
to the subring $\hxt\subs\hhxt$.
We therefore have $$A = \widehat{A}\cap\hxt^W = ann(e).$$
Just as in the rationalized case, we have $\K_G = i_W^*\circ\K_T\circ\rgt$, 
hence
$$\Im(\K_G)\cong i_W^*\left(\Im\K_T\circ\rgt\right)
\cong \frac{(\Im\K_T)^W}{ann(e)}.$$

Now consider the analogous diagram
$$
\xymatrix{
0 \ar[r] & B \ar[r]\ar[d] & \hxt \ar[r]^{i^*} \ar[d]
& \hso(\muh^{-1}(\xi,0)/T) \ar[d]\\
0 \ar[r] & \widehat{B} \ar[r] & \hhxt \ar[r]^{i^*} & \hso(\muh^{-1}(\xi,0)/T).}
$$
Since we have not assumed that $\muh^{-1}(\xi,0)/T$ is equivariantly formal,
we only know that the first two downward arrows are inclusions, and hence
can only conclude that $B$ is contained in the annihilator of $e'$.
Since $e'$ divides $e$, we have a series of natural surjections
$$\frac{\imkt^W}{ann(e)}\cong
\frac{\imkt^W}{A}\cong
\left(\frac{\Im\K_T}{B}\right)^W\to
\left(\frac{\Im\K_T}{ann(e')}\right)^W\to
\left(\frac{\Im\K_T}{ann(e)}\right)^W.$$
The composition of these maps is an isomorphism, hence so is each one.
\end{proofordinary}
\end{section}

\begin{section}{Quiver varieties}\label{quiver}
Let $Q$ be a quiver with vertex set $I$ and edge set $E\subs I\times I$,
where $(i,j)\in E$ means that $Q$ has an arrow pointing from $i$ to $j$.
We assume that
$Q$ is connected and has no oriented cycles.
Suppose given two collections of vector spaces $\{V_i\}$ and $\{W_i\}$,
each indexed by $I$, and consider the affine space
$$A = \bigoplus_{(i,j)\in E}\Hom(V_i,V_j)\oplus\bigoplus_{i\in I}\Hom(V_i,W_i).$$
The group $U(V) = \prod_{i\in I}U(V_i)$ acts on $A$ by conjugation, and
this action is hamiltonian.  Given an element
$$(B,J) = \bigoplus_{(i,j)\in E}B_{ij}\oplus\bigoplus_{i\in I}J_i$$
of $A$,
the $\mathfrak{u}(V_i)^*$ component of the moment map
is $$\mu_i(B,J) = J_i^{\dagger}J_i+\sum_{(i,j)\in E}B_{ij}^{\dagger}B_{ij},$$
where $\dagger$ denotes adjoint, and $\mathfrak{u}(V_i)^*$
is identified with with the set of hermitian matrices via
the trace form.
Given a generic central element $\xi\in\mathfrak{u}(V)^*$,
the K\"ahler quotient $A\mod_{\!\xi}U(V)$ parameterizes stable, 
framed representations
of $Q$ of fixed dimension \cite{N1}.  If $W_i=0$ for all $i$, then
the diagonal circle $U(1)$ in the center of $U(V)$ acts trivially, 
and we instead quotient by
$PU(V) = U(V)/U(1)$.

Consider the hyperk\"ahler quotient $$\M = T^*A\mmod_{\!(\xi,0)}U(V).$$
As in Example~\ref{standard}, $\M$ has a natural circle action induced from scalar
multiplication on the fibers of $T^*A$.  
We now show that $X=T^*A$ satisfies
the hypotheses of Theorems~\ref{main} and \ref{ordinary}.

\begin{proposition}
Let $T(V)\subs U(V)$ be a maximal torus, and let $pr:\mathfrak{u}(V)^*\to\mathfrak{t}(V)^*$
be the natural projection.
The moment maps $\mu=\displaystyle{\bigoplus}_{i\in I}\mu_i:A\to\mathfrak{u}(V)^*$ 
and $pr\circ\mu:A\to\mathfrak{t}(V)^*$ are each proper.
\end{proposition}

\begin{proof}
To show that $\mu$ and $pr\circ\mu$ is proper, it suffices to find an element
$t\in T(V)\subs U(V)$ such that the weights of the 
action of $t$ on $A$ are all strictly positive.
Let $\lambda=\{\lambda_i\mid i\in I\}$ be a collection of integers,
and let $t\in T(V)$ be the central element of $U(V)$ that acts
on $V_i$ with weight $\lambda_i$ for all $i$.
Then $t$ acts on $\Hom(V_i,V_j)$ with weight $\lambda_j-\lambda_i$,
and on $\Hom(V_i,W_i)$ with weight $-\lambda_i$.
Hence we have reduced the problem to showing that it is possible
to choose $\lambda$ such that $\lambda_i<0$ for all $i\in I$ and
$\lambda_i<\lambda_j$ for all $(i,j)\in E$.

We proceed by induction on the order of $I$.
Since $Q$ has no oriented cycles, there must exist a source $i\in I$;
a vertex such that for all $j\in I$, $(j,i)\notin E$.
Deleting $i$ gives a smaller (possibly disconnected) quiver with
no oriented cycles, and therefore we may choose 
$\left\{\lambda_j\mid j\in I\smallsetminus\{i\}
\right\}$ such that $\lambda_j<0$ for all $j\in I\smallsetminus\{i\}$
and $\lambda_j<\lambda_k$ for all $(j,k)\in E$.
We then choose $\lambda_i<\min\left\{\lambda_j\mid j\in I\smallsetminus\{i\}\right\}$,
and we are done.
\end{proof}

\vspace{-\baselineskip}
\begin{proposition}
The rationalized Kirwan map 
$\Kh_{U(V)}:\widehat{H}^*_{\so\times U(V)}(T^*A)\to\hhso(\M)$
is surjective.
\end{proposition}

\begin{proof}
Nakajima \cite[\S 7.3]{N2} shows that there exist cohomology classes 
$a_i, b_i$ in the image of $\Kh_{U(V)}$ 
such that $\dso = \sum\pi_1^*a_i\cdot\pi_2^*b_i$.
(Nakajima uses a slightly modified circle action, but
his proof is easily adapted to the circle action that we have defined.)
It follows from Proposition~\ref{basis} that 
the classes $\{b_i\}$ generate $\hhso(\M)$.
\end{proof}

\vspace{-\baselineskip}
\begin{remark}\label{generation} 
This Proposition shows that the assumptions of Theorems~\ref{integration},~\ref{main}, 
and~\ref{ordinary} are satisfied for Nakajima's quiver varieties. 
Thus integration in equivariant cohomology yields a
description of the rationalized $\so$-equivariant cohomology, and also of 
the image of the non-rationalized Kirwan map $\K_G$. Therefore 
if we know that $\K_G$ is surjective for a particular quiver variety, 
then we have a concrete description of the ($\so$-equivariant) 
cohomology ring of that quiver variety. 
It is known that $\K_G$ is surjective
for Hilbert schemes of $n$ points on an ALE space, 
so our theory applies and gives a description
of the cohomology ring of these quiver varieties. 
It would be interesting to compare
our result in this case with that of \cite{LS} and \cite{LQW}. 
More examples of quiver varieties with surjective Kirwan map
are given in Remark \ref{horn}.
\end{remark}

\vspace{-\baselineskip}
\begin{remark}
Another interesting application of Proposition~\ref{basis} is 
for the moduli space of Higgs bundles. 
It is an easy exercise to write down the cohomology class of the diagonal in ${\mathcal M}\times {\mathcal M}$ 
as an expression in the tautological classes for the equivariantly formal and circle compact 
moduli space ${\mathcal M}$ of stable rank $n$ and degree $1$ 
Higgs bundles  on a genus $g>1$ smooth projective algebraic curve $C$. Therefore Proposition~\ref{basis} implies that 
the rationalized $\so$-equivariant cohomology ring $\hhso({\mathcal M})$ is generated by tautological classes.
In fact the same result follows from the argument of \cite{HT1}. There ${\mathcal M}$ was embedded into a circle compact manifold
${\mathcal M}_\infty$, whose cohomology is the free algebra on the tautological classes. The argument in \cite{HT1} then goes by
showing that the embedding of the $\so$-fixed point set of ${\mathcal M}$ in that of ${\mathcal M}_\infty$ induces a surjection
on cohomology. This already implies that $\hhso({\mathcal M}_\infty)$ surjects onto $\hhso(\mathcal M)$. In \cite{HT1} it is shown
that in the rank $2$ case this embedding also implies the surjection on non-rationalized cohomology, and then
a companion paper \cite{HT2} describes the cohomology ring of ${\mathcal M}$ explicitly. However for higher rank
this part of the argument of \cite{HT1} breaks down. Later Markman \cite{Ma} used similar diagonal arguments on 
certain compactifications of $\mathcal M$ and
 Hironaka's celebrated theorem on desingularization of algebraic varieties to deduce that the cohomology ring of ${\mathcal M}$ 
is generated by tautological classes for all $n$. 
\end{remark}

\begin{example}\label{thaddeus}
Here we present an example of an embedding of circle compact manifolds, due to Thaddeus \cite{T}, 
where surjection on rationalized $\so$-equivariant cohomology does not imply surjection on $\so$-equivariant cohomology. 
Let $\so$ act on  $\P^1\times\P^1$ by 
$$((x:y),(u:v))\mapsto ((\lambda x :y),(u:v))$$ and on $\P^3$ by 
$$(z_1:z_2:z_3:z_4)\mapsto (\lambda z_1: \lambda z_2: z_3:z_4).$$ 
Then the Segr\'e embedding $i:\P^1\times\P^1\to \P^3$ given by
$$i\big((x:y),(u:v)\big)=(xu:xv:yu:yv)$$ 
is $\so$-equivariant, and clearly induces an isomorphism 
on the fixed point sets of the
$\so$ action. Therefore $i^*: \hhso(\P^3)\to \hhso(\P^1\times \P^1)$ is surjective, and in fact an isomorphism, however
$i^*:\hso(\P^3)\to \hso(\P^1\times \P^1)$ is only an injection and therefore cannot be surjective. 
\end{example}

\begin{section}{Hyperpolygon spaces}\label{polygon} 

We conclude by illustrating Theorem~\ref{ordinary} with
a computation of the equivariant cohomology ring of a hyperpolygon space.
Proposition~\ref{hp} first appeared in \cite{HP2}, and Corollary~\ref{konno}
in \cite{K2}, both obtained by geometric arguments completely
different from those used here.

A hyperpolygon space, introduced in \cite{K2},
is a quiver variety associated to the 
following quiver (Figure \ref{fig:quiver}),
with $V_0 = \C^2$, $V_i=\C^1$ for $i\in\{1,\ldots,n\}$, and $W_i=0$ for all $i$.
It is so named because, for 
$$\xi = \left(-\half\sum_{i=1}^n\xi_i;\hspace{2pt}\xi_1,\ldots,\xi_n\right)
\in\mathfrak{pu}(V)^*\subs\mathfrak{u}(2)^*\oplus\mathfrak{u}(1)^n,$$
the K\"ahler quotient $A\mod_{\!\xi}PU(V)\cong(\C^2)^n\mod_{\!\xi}PU(V)$ 
parameterizes $n$-sided polygons in $\R^3$
with edge lengths $\xi_1,\ldots,\xi_n$, up to rotation \cite{HK}.
\begin{figure}[h]
\centerline{\epsfig{figure=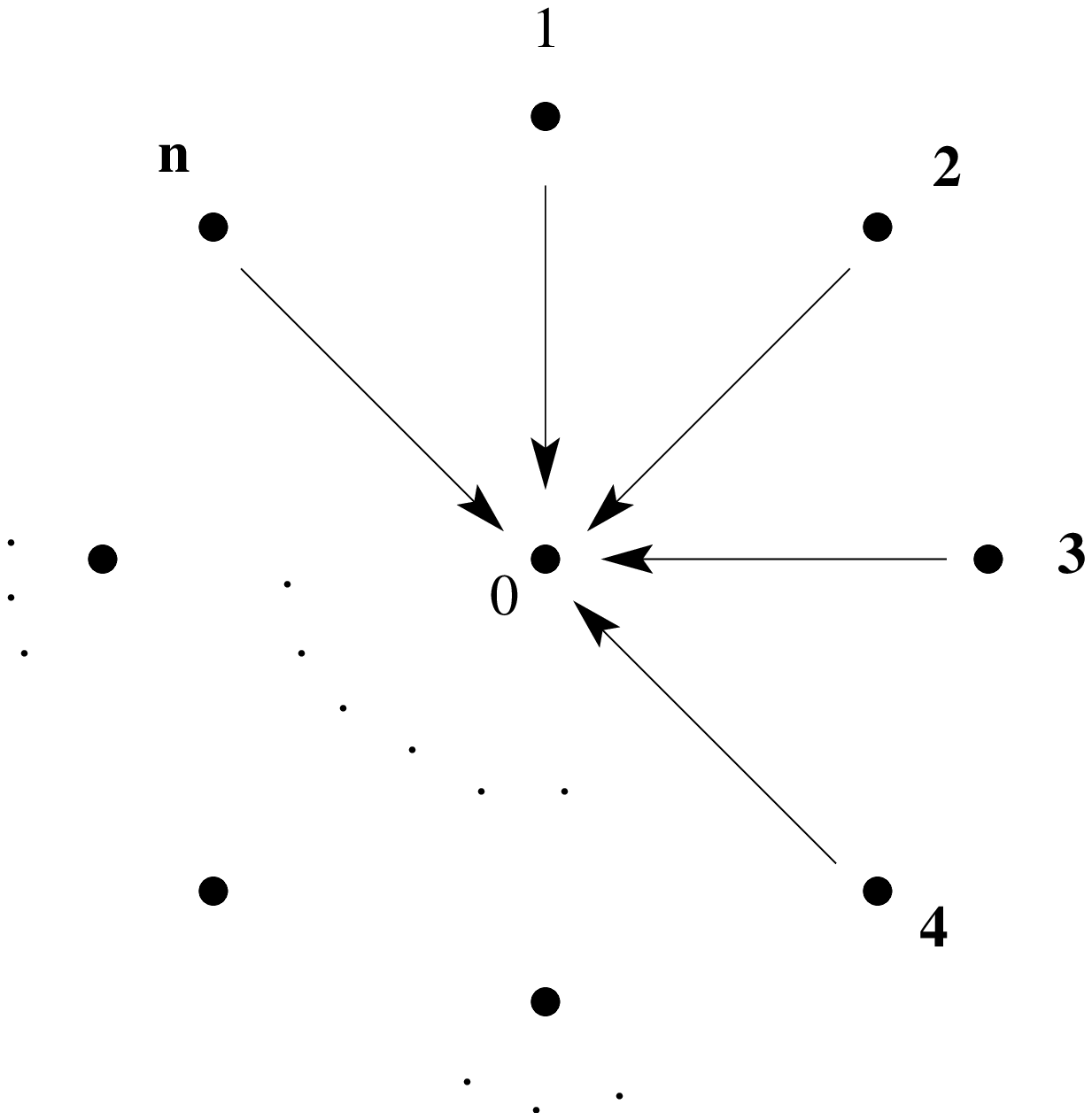,height=4cm}}
\caption{The quiver for a hyperpolygon space.}\label{fig:quiver}
\end{figure}

We will simplify our computations by dividing first by
the torus $\prod_{i=1}^nU(V_i)$.
We have
\begin{eqnarray*}
\M &=& \left(T^*\C^2\right)^n\bigmmod PU(V)\\
&\cong&\left(\left(T^*\C^2\right)^n\bigmmod\prod_{i=1}^nU(V_i)\right)\bigmmod SU(2)\\
&\cong& \prod_{i=1}^nT^*\C P^1\bigmmod SU(2),
\end{eqnarray*}
where the action of $SU(2)$ on each copy of $T^*\C P^1$
is induced by the rotation action on $\C P^1\cong S^2$.

\begin{proposition}\label{hpsurj}
The non-rationalized Kirwan map
$\K_{U(V)}:H^*_{\so\times U(V)}(T^*\C^{2n})\to\hso(\M)$ is surjective.
\end{proposition}

\begin{proof}
The map $\K_{U(V)}$ factors as a composition
$$H^*_{\so\times U(V)}(T^*\C^{2n})
\to H^*_{\so\times SU(2)}\left(\prod_{i=1}^nT^*\C P^1\right)
\overset{\K_{SU(2)}}\longrightarrow \hso(\M),$$
where the first map is the Kirwan map for a toric hyperk\"ahler variety,
and therefore surjective by \cite{HP1}.  Hence it suffices to show that
$\K_{SU(2)}$ is surjective.

The level set $\muc^{-1}(0)$ for the action
of $SU(2)$ on $\prod_{i=1}^nT^*\C P^1$ is a subbundle of the cotangent bundle,
given by requiring the $n$ cotangent vectors to add to zero after being restricted
to the diagonal $\C P^1$.  In particular this set is smooth, 
and its $\so\times SU(2)$-equivariant cohomology ring is canonically
isomorphic to that of $\prod_{i=1}^nT^*\C P^1$.
Hence $\K_{SU(2)}$ factors as
$$H^*_{\so\times SU(2)}\left(\prod_{i=1}^nT^*\C P^1\right)
\cong H^*_{\so\times SU(2)}\big(\muc^{-1}(0)\big)\to
\hso\big(\muc^{-1}(0)\mod SU(2)\big)
\cong \hso(\M),$$ where the map in the middle is the K\"ahler
Kirwan map.
Surjectivity of this map follows from the following more general lemma,
applied to the manifold $\muc^{-1}(0)$.

\begin{lemma}\label{cut}
Suppose given a hamiltonian action of $\so\times G$ on a symplectic
manifold $X$, such that the $\so$ component of the moment map is proper
and bounded below with finitely many critical values.  
Then the K\"ahler Kirwan map
$\K:H^*_{\so\times G}(X)\to\hso(X\mod G)$ is surjective.
\end{lemma}

\begin{proof}
Extend the action of $S^1$ to an action on $X\times\C$ by letting
$S^1$ act only on the left-hand factor.  On the other hand,
consider a second copy of the circle, which we will call $\mathbb T$ to avoid confusion, 
acting diagonally on $X\times\C$.
Choose $r\in \operatorname{Lie}(\mathbb T)^*\cong\R$ greater than
the largest critical value of the $\mathbb T$-moment map, and consider
the space $$Cut(X\mod G):=\left(X\times\C\right)\mod_{\! r} \mathbb T\times G
\cong\big((X\mod G)\times\C\big)\mod_{\! r} \mathbb T.$$
This space, which is called the {\em symplectic cut} of $X\mod G$ \cite{L},
is an $S^1$-equivariant (orbifold) compactification of $X\mod G$.
We then have a commutative diagram
$$\begin{CD}
H^*_{\so\times G\times \mathbb T}(X\times\C)        @>>> H^*_{\so\times G}(X)\\
@VVV @VV\K V\\
\hso(Cut(X\mod G))        @>>> \hso(X\mod G).
\end{CD}$$
The vertical map on the left is surjective because the $G\times \mathbb T$ moment map is proper,
and the map on the bottom is surjective because the long
exact sequence in cohomology for $X\mod G\subs Cut(X\mod G)$
splits naturally, hence $\K$ is surjective as well.
\end{proof}

\noindent By applying Lemma~\ref{cut} to $X = \muc^{-1}(0)$,
this completes the proof of Proposition~\ref{hpsurj}.
\end{proof}

\vspace{-\baselineskip}
\begin{remark}\label{horn}
The argument in Proposition \ref{hpsurj} generalizes immediately
to show that the hyperk\"ahler Kirwan map for the quotient
$$\left(\prod_{i=1}^nT^*Flag(\C^k)\right)\bigmmod SU(k)$$
is surjective.  This is itself a quiver variety, and like the
hyperpolygon space, it has a moduli theoretic interpretation.
The K\"ahler quotient
$$\left(\prod_{i=1}^nFlag(\C^k)\right)\bigmod SU(k)$$
is isomorphic to the space of $n$-tuples of $k\times k$
hermitian matrices with fixed eigenvalues adding to zero, modulo 
conjugation.
This space has been studied by many authors.  The classical
problem, due to Horn, of determining the values of 
the moment map for which it is
nonempty, has only recently been solved \cite{KT}.  For a survey, see \cite{Fu}.
\end{remark}

To compute the kernel of the hyperk\"ahler Kirwan map for the
hyperpolygon space, we first need to study the abelian quotient
$$\NN := \prod_{i=1}^nT^*\C P^1\bigmmod T,$$ where $T\cong U(1)\subs SU(2)$ is a
maximal torus. 
The space $\prod_{i=1}^nT^*\C P^1$ is a toric hyperk\"ahler manifold
\cite{BD}, given by an arrangement of $2n$ hyperplanes in $\R^n$,
where the $(2i-1)^{\text{st}}$ and $(2i)^{\text{th}}$ hyperplanes
are given by the equations $x_i = \pm\xi_i$.
Taking the hyperk\"ahler quotient by $T$ corresponds on the level of arrangements
to restricting this arrangement to the hyperplane
$\{x\in\R^n\mid\sum x_i = 0\}$.

Call a subset $S\subs\{1,\ldots,n\}$ {\em short} if $\Sd_{i\in S}\xi_i<\Sd_{j\in S^c}\xi_j$.
Requiring that $\xi$ is a regular value of the hyperk\"ahler moment
map is equivalent to requiring that for every $S\subs\{1,\ldots,n\}$,
either $S$ or $S^c$ is short \cite{K2}.
Applying \cite[4.5]{HP1}, we have
$$\hso\left(\NN\right)\cong
\Q[a_1,b_1,\ldots,a_n,b_n,\a,x]\Big/
\Big\la a_i-b_i-\a, \hs a_ib_i\hs\Big{|}\hs i\leq n\Big\ra+
\Big\la A_S, B_S\hs\Big{|}\hs S\text{   short}\Big\ra,
$$
where 
\begin{equation*}\label{asbs}
A_S= \prod_{i\in S}(x-a_i)\prod_{j\in S^c}b_j\hs\hs
\text{   and   }\hs\hs B_S = \prod_{i\in S}(x-b_i)\prod_{j\in S^c}a_j.
\end{equation*}
Here $\a$ is the image in $\hso\left(\NN\right)$
of the unique positive root of $SU(2)$.
The Weyl group $W$ of $SU(2)$, isomorphic to $\Z/2\Z$, acts on this ring by fixing $x$
and switching $a_i$ and $b_i$ for all $i$.
Let $c_i = a_i+b_i$, and let $C_S=A_S+B_S$.
Let $$P=\Q[c_1,\ldots,c_n,\a,x]\Big/
\Big\la c_i^2-\a^2\hs\Big{|}\hs i\leq n\Big\ra$$ and
$$Q=P^W=\Q[c_1,\ldots,c_n,\a^2,x]\Big/
\Big\la c_i^2-\a^2\hs\Big{|}\hs i\leq n\Big\ra.$$
Let $$\I=\Big\la A_S, B_S\hs\Big{|}\hs S\text{   short}\Big\ra\subs P
\hspace{15pt}\text{and}\hspace{15pt} 
\J=\Big\la C_S\hs\Big{|}\hs S\text{   short}\Big\ra\subs Q,$$
so that $$\hso(\NN)\cong P/\I
\hspace{15pt}\text{and}\hspace{15pt} 
\hso(\NN)^W\cong Q/\J.$$
Note that all odd powers of $\a$ in the expression for $C_S=A_S+B_S$ cancel out.

Then by Theorem~\ref{ordinary} and Remark~\ref{notsobad},
\begin{eqnarray*}
\hso(\M)&\cong& \frac{\hso\left(\NN\right)^W}{ann(e)}
\cong \frac{Q}{(e:\J)},
\end{eqnarray*}
where $e=\a^2(x^2-\a^2)$, and $(e:\J)$ is the ideal of elements
of $Q$ whose product with $e$ lies in $\J$.

If $S$ is a nonempty short subset, let
$m_S$ be the smallest element of $S$, $n_S$ the smallest
element of $S^c$, and $$D_S = \prod_{m_S\neq i\in S}(c_i-x)
\cdot\prod_{n_S\neq j\in S^c}(c_{n_S}+c_j)\in Q.$$

\begin{proposition}\label{hp}
The equivariant cohomology ring $\hso(\M)$ is isomorphic to\footnote{The class
denoted by $c_i$ in \cite{HP2} differs from our $c_i$ by a sign, hence
to recover the presentation of \cite{HP2} we must replace $c_i-x$
with $c_i+x$ in the expression for $D_S$.}
$$Q\big/\big\la
D_S\mid \emptyset\neq S\text{   short}\big\ra.$$
\end{proposition}

\begin{proof}
We begin by proving that $e\cdot D_S\in\mathcal{J}$ for all nonempty
short subsets $S\subs\{1,\ldots,n\}$.  We will in fact prove the slightly
stronger statement
$$e\cdot D_S\in\Big\la C_T\hs\Big{|}\hs T\subs S\text{   short}\Big\ra
\subs\mathcal{J},$$
proceeding by induction on $|S|$.
We will assume, without loss of generality, that $n\in S$.
The base case occurs when $S=\{n\}$, and in this case
we observe that
$$e\cdot D_S = 2^{n-3}\cdot(x+c_n)\cdot\big((2x-c_n)\cdot C_{\emptyset} - c_n\cdot C_S\big).$$
We now proceed to the inductive step, assuming that the proposition
is proved for all short subsets of size less than $|S|$, and all values of $n$.
For all $T\subs S\smallsetminus\{n\}$, we have
$$\half\big(C_T-C_{T\cup\{n\}}\big)=(c_n-x)\cdot C'_T,$$
where $C'_T$ is the polynomial in the variables $\{c_1,\ldots,c_{n-1},\a^2\}$
corresponding to the short subset $T\subs\{1,\ldots,n-1\}$.
Since $S\smallsetminus\{n\}$ is a short subset of $\{1,\ldots,n-1\}$
of size strictly smaller than $S$,
our inductive hypothesis tells us that $e\cdot D_S/(c_n-x)$
can be written as a linear combination of polynomials $C'_T$, where the coefficients
are quadratic polynomials in $\{c_1,\ldots,c_{n-1},\a^2\}$.
Replacing $C'_T$ with $\half\left(C_T-C_{T\cup\{n\}}\right)=(c_n-x)\cdot C'_T$, 
we have expressed $e\cdot D_S$ in terms of the appropriate polynomials.
This completes the induction.

Suppose that $F\in Q$ is an element of degree less than $n-2$
such that $e\cdot F\in\J$.  By the second isomorphism
of Theorem~\ref{ordinary}, this implies that $e'\cdot F\in\I\subs P$,
where $e'=\a(x^2-\a^2)$.
Consider the quotient ring $R$ of $P$ obtained by setting
$a_i^2=b_i^2=x=0$ for all $i$.  (Recall that $a_i=\half(c_i+\a)$
and $b_i=\half(c_i-\a)$.)  Then element $e'$ maps to zero in $R$,
while the generators $\{A_S, B_S\}$ of $\I$ descend to a basis
for the $n^{\text{th}}$ degree part of $R$.  This means that we
must have $e'\cdot F=0\in P$.
Using the additive basis for $P$ consisting of monomials
that are squarefree in the variables $c_1,\ldots,c_n$,
it is easy to check that $e'$ is not a zero divisor in $P$,
and therefore that $F=0$.

Finally, we must show that $\{D_S\mid \emptyset\neq S\text{  short}\}$
generates all elements of $(e:\J)$ of degree at least $n-2$.
Let $F$  be an element of minimal degree $k\geq n-2$ that is in
$(e:\J)$ but not $\la D_S\mid \emptyset\neq S\text{  short}\ra$.
In the proof of \cite[3.2]{HP2} it is shown that
$\{D_S\mid \emptyset\neq S\text{  short}\}$ descends to a basis for the degree $n-2$
part of the quotient ring $Q/\la x\ra$, hence $F$
differs from an element of 
$\la D_S\mid \emptyset\neq S\text{  short}\ra$ by $x\cdot F'$ for some $F'$
of degree $k-1$.  By equivariant formality of $\hso(\M)$,
$$x\cdot F' = F\in (e:\J)\impl F'\in (e:\J),$$
which contradicts the minimality of $k=\operatorname{deg}F$.
Hence $\la D_S\mid \emptyset\neq S\text{  short}\ra=(e:\J)$, and the proposition
is proved.
\end{proof}
\end{section}

\vspace{-\baselineskip}
\begin{corollary}\label{konno}
The ordinary cohomology ring $H^*(\M)$ is isomorphic to 
$$\Q[c_1,\ldots,c_n]\Big/\big\la c_i^2-c_j^2\mid i,j\leq n\big\ra
+\la\text{all monomials of degree $n-2$}\ra.$$
\end{corollary}

\begin{proof}
This follows from the fact that $H^*(M)\cong \hso(M)/\la x\ra$
for any equivariantly formal space $M$, and the observation in \cite{HP2}
that $\{D_S\mid \emptyset\neq S\text{  short}\}$ descends to a basis for the degree $n-2$
part of $Q/\la x\ra$.
\end{proof}
\end{section}

\footnotesize{
}
\end{spacing}
\end{document}